\documentclass[aap,MSNbibl,dvips]{arximspdf}
\usepackage{graphicx}
\doi{10.1214/13-AAP948} 
\volume{24}
\issue{3}
\pubyear{2014}
\firstpage{1269}
\lastpage{1297}

\makeatletter
\newcommand{\eqref}[1]{(\ref{#1})}
\newtheorem{theorem}{Theorem} 
\newtheorem{lemma}[theorem]{Lemma}

\newtheorem{proposition}[theorem]{Proposition}
\newproclaim{remark}[theorem]{Remark}
\newproclaim{example}[theorem]{Example}
\newproclaim{definition}[theorem]{Definition}
\newproclaim{notation}[theorem]{Notation}
\newproclaim{hypothesis}[theorem]{Hypothesis}
\newcommand{\bF}{\mathbb{F}}
\newcommand{\bB}{\mathbb{B}}
\newcommand{\bL}{\mathbb{L}}
\newcommand{\bN}{\mathbb{N}}
\newcommand{\bR}{\mathbb{R}}
\newcommand{\bV}{\mathbb{V}}

\newcommand{\cB}{\mathcal{B}}
\newcommand{\cF}{\mathcal{F}}
\newcommand{\cG}{\mathcal{G}}
\newcommand{\cL}{\mathcal{L}}
\newcommand{\cT}{\mathcal{T}}

\newcommand{\unif}{\operatorname{unif}}
\newcommand{\Bin}{\operatorname{Bin}}
\newcommand{\Beta}{\operatorname{Beta}}
\newcommand{\IPL}{\operatorname{IPL}}
\newcommand{\WI}{\operatorname{WI}}
\newcommand{\Sil}{\operatorname{Sil}}
\newcommand{\mSil}{\operatorname{mSil}}
\makeatother

\begin{document}
\begin{frontmatter}

\title{Search trees: Metric aspects and strong~limit~theorems}
\runtitle{Metric search trees}

\begin{aug}
\author[A]{\fnms{Rudolf} \snm{Gr\"ubel}\corref{}\ead[label=e1]{rgrubel@stochastik.uni-hannover.de}}
\runauthor{R. Gr\"ubel}
\affiliation{Leibniz Universit\"at Hannover}
\address[A]{Institut f\"ur Mathematische Stochastik\\
Leibniz Universit\"at Hannover\\
Postfach 6009\\
30060 Hannover\\
Germany\\
\printead{e1}} 
\end{aug}

\received{\smonth{2} \syear{2012}}
\revised{\smonth{6} \syear{2013}}


\begin{abstract}
We consider random binary trees that appear as the output of certain
standard algorithms for sorting and searching if
the input is random. We introduce the subtree size metric on search
trees and show that the resulting metric spaces
converge with probability~1. This is then used to obtain almost sure
convergence for various tree functionals, together
with representations of the respective limit random variables as
functions of the limit tree.
\end{abstract}

%

\pdfkeywords{60B99, 60J10, 68Q25, 05C05, Doob-Martin compactification, metric trees, path length, silhouette, subtree size metric, vector-valued martingales, Wiener index}
\end{frontmatter}

\section{Introduction}\label{sec:intro}
A sequential algorithm transforms an input sequence $t_1,t_2,\ldots$
into an output sequence $x_1,x_2,\ldots$ where, for
all $n\in\bN$, $x_{n+1}$ depends on $x_n$ and $t_{n+1}$ only.
Typically, the output variables are elements of some
combinatorial family $\bF$, each $x\in\bF$ has a size parameter
$\phi
(x)\in\bN$ and $x_n$ is an element of the set
$\bF_n:=\{x\in\bF\dvtx \phi(x)=n\}$ of objects of size $n$. In the
probabilistic analysis of such algorithms, one
starts with a stochastic model for the input sequence and is interested
in certain aspects of the output sequence. The
standard input model assumes that the $t_i$'s are the values of a
sequence $\eta_1,\eta_2,\ldots$ of independent and
identically distributed random variables. For random input of this
type, the output sequence then is the path of a Markov
chain $X=(X_n)_{n\in\bN}$ that is adapted to the family $\bF$ in the
sense that
%
\begin{equation}
\label{eq:time-space} P(X_n\in\bF_n) = 1 \qquad\mbox{for all } n\in
\bN.
\end{equation}
Clearly, $X$ is highly transient---no state can be visited twice.

The special case we are interested in, and which we will use to
demonstrate an approach that is generally applicable in
the situation described above, is that of binary search trees and two
standard algorithms, known by their acronyms BST
(binary search tree) and DST (digital search tree). These are discussed
in detail in the many excellent texts in this area,
for example in \cite{Knuth3,Mahmoud1} and \cite{Drmota09}.
Various functionals of the search trees, such as the
height \cite{Dev86}, the path length \cite{RegnierQS,RoeQS}, the node depth
profile \cite{J-H2001,CDJ2001,CR2004,CKMR2005,FHN2006,DJN2008},
the subtree size profile \cite{Fu2008,DeGr3}, the Wiener index \cite{Nein02} and the
silhouette \cite{GrSilh} have been studied, with methods spanning the
wide range from generatingfunctionology to
martingale methods to contraction arguments on metric spaces of
probability distributions (neither of these lists is
complete). Many of the results are asymptotic in nature, where the
convergence obtained as $n\to\infty$ may refer to the
distributions or to the random variables themselves. As far as strong
limit theorems are concerned, a significant step
toward a unifying approach was made in the recent paper \cite{EGW},
where methods from discrete potential theory were
used to obtain limit results on the level of the combinatorial
structures themselves: In a suitable extension
of the state space $\bF$, the random variables $X_n$ converge almost
surely as $n\to\infty$, and the limit generates the
tail $\sigma$-field of the Markov chain. The results in \cite{EGW}
cover a wide variety of structures; search trees are
a special case. It should also be mentioned here that the use of
boundary theory has a venerable tradition in connection
with random walks; see \cite{KaimVersh} and~\cite{Woess1}.

Our aims in the present paper are the following. First, we use the
algorithmic background for a direct proof of the
convergence of the BST variables $X_n$, as $n\to\infty$, to a limit
object $X_\infty$, and we obtain a representation
of $X_\infty$ in terms of the input sequence $(\eta_i)_{i\in\bN}$.
Second, we introduce the subtree size metric on
finite binary trees. This leads to a reinterpretation of the above
convergence in terms of metric trees. We also
introduce a family of weighted variants of this metric, with parameter
$\rho\ge1$, and then identify the critical
value $\rho_0$ with the property that the metric trees converge for
$\rho<\rho_0$ and do not converge if $\rho>
\rho_0$. The value $\rho_0$ turns out to also be the threshold for
compactness of the limit tree. Third, we use
convergence at the tree level to (re)obtain strong limit theorems for
three tree functionals---the path length, the
Wiener index and a metric version of the silhouette.

These topics are treated in the next three sections, where each has its
own introductory remarks.

\section{Binary search trees}%
\label{sec:BSTDoob--Martin}

We first introduce some notation, mostly specific to binary trees, then
discuss the two search algorithms and the associated
Markov chains and finally recall the results from \cite{EGW} related to
these structures, including an alternative proof
of the main limit theorem.

\subsection{Some notation}\label{subsec:notation}
We write $\cL(X)$ for the distribution of a random variable $X$ and
$\cL
(X|Y=k)$, $\cL(X|Y)$, $\cL(X|\cF)$ for the
various versions of the conditional distribution of $X$ given (the
value of) a random variable $Y$ or a $\sigma$-field
$\cF$. Further, $\delta_c$ is the one-point mass at $c$, $1_A$ is the
indicator function of the set $A$ [so that
$1_A(c)=\delta_c(A)$], $\Bin(n,p)$ denotes the binomial distribution
with parameters $n\in\bN$ and $p\in(0,1)$,
$\Beta(\alpha,\beta)$ is the beta distribution with parameters
$\alpha
,\beta>0$ and $\unif(0,1)=\Beta(1,1)$ is the
uniform distribution on the unit interval. We also write $\unif(M)=(\#
M)^{-1}\sum_{c\in M} \delta_c$ for the uniform
distribution on a finite set $M$.

With $\bN_0=\{0,1,2,\ldots\}$ let
\[
\bV_k:= \{0,1\}^k,\qquad \bV:=\bigsqcup
_{k\in\bN_0}\bV_k,\qquad \partial\bV:= \{0,1\}^\infty
\]
be the set of 0--1 sequences of length $k$, $k\in\bN_0$, the set of all
finite 0--1 sequences and the set of all infinite
0--1 sequences, respectively. The set $\bV_0$ has $\varnothing$, the
``empty sequence,'' as its only element, and $|u|$ is the
length of $u\in\bV$, that is, $|u|=k$ if $u\in\bV_k$. For each node
$u=(u_1,\ldots,u_k)\in\bV$ we use
\begin{eqnarray*}
u0&:=&(u_1,\ldots,u_k,0),
\\
u1&:=&(u_1,\ldots,u_k,1),
\\
\bar u&:=&(u_1,\ldots,u_{k-1})\qquad \mbox{if } k\ge1,
\end{eqnarray*}
to denote its left and right direct descendant (child) and its direct
ancestor (parent). We write $u\le v$ for
$u=(u_1,\ldots,u_k)\in\bV$, $v=(v_1,\ldots,v_l)\in\bV$ if $k\le
l$ and
$u_j=v_j$ for $j=1,\ldots,k$, that is, if $u$ is a
prefix of $v$; the extension to $v\in\partial\bV$ is obvious. The
prefix order is a partial order only, but there exists
a unique minimum $u\wedge v$ to any two nodes $u,v\in\bV$, their last
common ancestor; again, this can be extended to
elements of $\partial\bV$. Another ordering on $\bV$ can be obtained
via the function $\beta\dvtx\bV\to[0,1]$,
%
\begin{equation}
\label{eq:defbeta} \beta(u):= \frac{1}{2} + \sum_{j=1}^k
\frac
{2u_j-1}{2^{j+1}}, \qquad u\in\bV.
\end{equation}
This will be useful in various proofs, and also in connection with
illustrations.

By a \emph{binary tree} we mean a subset $x$ of the set $\bV$ of nodes
that is prefix stable in the sense that $u\in x$ and
$v\le u$ implies that $v\in x$. Informally, we regard the components
$u_1,\ldots,u_k$ of $u$ as a routing instruction
leading to the vertex $u$, where~0 means a move to the left, 1 a move
to the right and the empty sequence is the root
node. The edges of the tree $x$ are the pairs $(\bar u,u)$, $u\in x,
u\neq\varnothing$. A node is external to
a tree if it is not one of its elements, but its direct ancestor is; we
write $\partial x:=\{u\in\bV\dvtx \bar u\in x,
u\notin x\}$ for the set of external nodes of $x$. Finally,
%
\begin{equation}
\label{eq:sts} \sigma(x,u):= \#\{v\in x\dvtx u\le v\}
\end{equation}
is the size of the subtree of $x$ rooted at $u$ (or the number of
descendants of $u$ in $x$, including $u$).

Let $\bB$ denote the (countable) set of finite binary trees, $\bB
_n:=\{
x\in\bB\dvtx \#x=n\}$ those of size (number of
nodes) $n$. The single element of $\bB_1$ is $\{\varnothing\}$, the tree
that consists of the root node only.

\subsection{Search algorithms and Markov chains}\label{subsec:AlgMC}
Let $(t_i)_{i\in\bN}$ be a sequence of pairwise distinct real numbers.
The BST (binary search tree) algorithm stores
these sequentially into labeled binary trees $(x_n,L_n)$, $n\in\bN$,
with $x_n\in\bB_n$ and
$L_n\dvtx x_n\to\{t_1,\ldots,t_n\}$. For $n=1$ we have $x_1=\{
\varnothing\}$
and $L_1(\varnothing)=t_1$. Given $(x_n,L_n)$, we
construct $(x_{n+1},L_{n+1})$ as follows: Starting at the root node we
compare the next input value $t_{n+1}$ to the
value $L_n(u)$ attached to the node $u$ under consideration, and move
to $u0$ if $t_{n+1}<L_n(u)$ and to $u1$
otherwise, until an ``empty'' node $u$ (necessarily an external node of
$x_n$) is found. Then $x_{n+1}:=x_n\cup\{u\}$ and
$L_{n+1}(u):=t_{n+1}$, $L_{n+1}(v):=L_n(v)$ for all $v\in x_n$.

Now let $(\eta_i)_{i\in\bN}$ be a sequence of independent random
variables with $\cL(\eta_i)=\unif(0,1)$ for all
$i\in\bN$, and let $X_n$ be the random binary tree associated with the
first $n$ of these. By construction, the label
functions $L_n$ are monotone with respect to the $\beta$-order of the
tree nodes, that is, with $\beta$ as
in \eqref{eq:defbeta},
%
\begin{equation}
\label{eq:orderpropBST} \beta(u)\le\beta(v) \Rightarrow L_n(u)\le
L_n(v)\qquad \mbox{for all } n \mbox{ with }\{u,v\}\subset
X_n.
\end{equation}
In particular, if we number the external nodes of $X_n$ from the left
to the right, then the number of the node that
receives $\eta_{n+1}$ is the rank of this value among $\{\eta
_1,\ldots
,\eta_n\}$, hence uniformly distributed on
$\{1,\ldots,n+1\}$. This shows that the (deterministic) BST algorithm,
when applied to the (random) input
$(\eta_i)_{i\in\bN}$, results in a Markov chain $(X_n)_{n\in\bN}$ with
state space $\bB$, start at $X_1\equiv
\{\varnothing\}$ and transition probabilities
%
\begin{equation}
\label{eq:transBST} Q \bigl(x,x\cup\{u\} \bigr) = %
\cases{1/(1+\# x), &\quad $\mbox{if }u\in\partial x$,\vspace*{2pt}
\cr
0, & \quad $\mbox{otherwise}$.} %
\end{equation}
In words: We obtain $X_{n+1}$ by choosing one of the $n+1$ external nodes
of $X_n$ uniformly at random and joining it to the tree. We refer to
this construction as the
\emph{BST chain}.

For the DST (digital search tree) algorithm, the input values are
infinite 0--1 sequences, that is, elements of $\partial
\bV$. Given $t_1,t_2,\ldots\in\partial\bV$ we again obtain a sequence
$x_1,x_2,\ldots$ of labeled binary trees, but now
we use the components $t_{n+1,k}$, $k\in\bN$, of the next input value
$t_{n+1}$ as a routing instruction through $x_n$,
moving to $u0$ from an occupied node $u\in\bV_k$ if $t_{n+1,k+1}=0$ and
to $u1$ otherwise. As in the BST case we assume
that the $t_i$'s are the values of a sequence of independent and
identically distributed random variables $\eta_i$,
where the distribution of the $\eta_i$'s is now a probability measure
$\mu$ on the measurable space
$(\partial\bV,\cB(\partial\bV))$, with $\cB(\partial\bV)$ the
$\sigma
$-field generated by the projections on the
sequence elements, $\partial\bV\ni t=(t_k)_{k\in\bN}\mapsto t_i$,
$i\in
\bN$. This $\sigma$-field is also generated by
the sets
%
\begin{equation}
\label{eq:defA} A_u:= \{v\in\partial\bV\dvtx v\ge u \}, \qquad u\in\bV.
\end{equation}
It is easy to check that the intersection of two such sets is either
empty or again of this form. This implies that
$\mu$ is completely specified by its values $\mu(A_u)$, $u\in\bV$, and
the DST analogue of \eqref{eq:transBST} then is
%
\begin{equation}
\label{eq:transDST} Q \bigl(x,x\cup\{u\} \bigr) = %
\cases{\mu(A_u), &\quad $\mbox{if } u\in\partial x$,\vspace*{2pt}
\cr
0, & \quad $\mbox{otherwise}$.} %
\end{equation}
By the \emph{DST chain with driving distribution $\mu$} we mean a
Markov chain $(X_n)_{n\in\bN}$ with state space $\bB$,
start at $\{\varnothing\}$ and transition mechanism given by \eqref
{eq:transDST}.

\subsection{Doob--Martin compactification}\label{subsec:DM}
We refer the reader to Doob's seminal paper \cite{Doob1959} and to the
recent textbook \cite{Woess2} for the main
results of, background on and further references for the boundary
theory for transient Markov chains. For the BST chain
the Doob--Martin\vspace*{1pt} compactification has recently been obtained in \cite{EGW}: It can be described as the closure $\bar\bB$
of the embedding of $\bB$ into the compact space $[0,1]^{\bV}$, endowed
with pointwise convergence, that is given by the
standardized \emph{subtree size functional}
\[
\bB\ni x \mapsto \biggl(\bV\ni u\mapsto\frac{\sigma(x,u)}{\# x} \biggr)
\]
with $\sigma$ as defined in \eqref{eq:sts}. Further, the elements of
the boundary $\partial\bB$ may be represented by
probability measures $\mu$ on $(\partial\bV,\cB(\partial\bV))$, with
convergence $x_n\to\mu$ of a sequence
$(x_n)_{n\in\bN}$ in $\bB$ meaning that
\[
\mu(A_u) = \lim_{n\to\infty} \frac{\sigma(x_n,u)}{\# x_n} \qquad\mbox{for all } u\in\bV,
\]
and $\mu_n(A_u)\to\mu(A_u)$ for all $u\in\bV$ if we have a sequence
$(\mu_n)_{n\in\bN}$ of elements of $\partial\bB$
instead.

The general theory implies that $X_n$ converges almost surely to a
limit $X_\infty$ with values in $\partial\bB$;
\cite{EGW} also contains a description of $\cL(X_\infty)$. The proof
given there does not make use of the algorithmic
background, but takes the transition mechanism \eqref{eq:transBST} as
its starting point. We now show that this
background leads to a direct proof of $X_n\to X_\infty$, and to a
representation of $X_\infty$ in terms of the input
sequence.

We need some more notation. On $\bV$ we define a metric $d_\bV$ by
%
\begin{equation}
\label{eq:defdV} d_\bV(u,v):= 2^{-|u\wedge v|}-\tfrac{1}{2}
\bigl(2^{-|u|}+2^{-|v|} \bigr),\qquad u,v\in\bV.
\end{equation}
On $\bV$ itself this gives the discrete topology, and the completion of
$\bV$ with respect to $d_\bV$ leads to $\bar
\bV:=\bV\cup\partial\bV$, a compact and separable metric space.
This is
also the ends compactification if we regard
$\bV$ as the complete rooted binary tree. We extend the $A_u$'s to
$\bar\bV$ by
\[
\bar A_u:= \{ v\in\bar\bV\dvtx v\ge u\}, \qquad u\in\bV.
\]
Because of
\[
\bar A_u:= \bigl\{ v\in\bar\bV\dvtx d_\bV(u,v)<2^{-|u|}
\bigr\} = \bigl\{ v\in\bar\bV\dvtx d_\bV(u,v)\le2^{-|u|-1}
\bigr\}
\]
these sets are open and closed. Further,
\[
\{u\} = \bar A_u\setminus(\bar A_{u0}\cup\bar
A_{u1}), \qquad \bar A_u\cap\bar A_v= %
\cases{ \bar A_u, &\quad $\mbox{if } u\le v$,\vspace*{2pt}
\cr
\bar
A_v, &\quad $\mbox{if } u\ge v$,\vspace*{2pt}
\cr
\varnothing, &\quad $
\mbox{otherwise,}$} %
\]
hence $\{\bar A_u\dvtx u\in\bV\}$ is a $\pi$-system that generates
$\cB
(\bar\bV)$. Together these facts imply that
weak convergence of probability measures $\mu_n$ to a probability
measure $\mu$ on $(\bar\bV,\cB(\bar\bV))$
is equivalent to
%
\begin{equation}
\label{eq:wconvT} \lim_{n\to\infty} \mu_n(\bar
A_u) = \mu(\bar A_u)\qquad \mbox{for all } u\in\bV.
\end{equation}
In view of
\[
\frac{1}{n} \sigma(X_n,u) = \unif(X_n) (\bar
A_u)
\]
and $X_\infty(\bV)=0$ convergence in the Doob--Martin topology is
therefore equivalent to the weak convergence of
probability measures on the metric space $(\bar\bV,d_\bV)$ if we
represent finite subsets $M$ of $\bV$ by the uniform
distribution $\unif(M)$ on $(\bar\bV,\cB(\bar\bV))$.

Moreover, any sequence $(\mu_n)_{n\in\bN}$ of probability measures on
$(\bar\bV,\cB(\bar\bV))$ is tight, as $\bar\bV$
is compact, and therefore has a limit point by Prohorov's theorem
\cite{Bill68}, page~37. If $(\mu_n(\bar A_u))_{n\in\bN}$
is a convergent sequence for each $u\in\bV$, then there is only one
such limit point, which means that $\mu_n$
converges weakly to some probability measure $\mu$ and that \eqref
{eq:wconvT} holds.
Finally, let
%
\begin{equation}
\label{eq:deftau} \tau(u):=\inf\{n\in\bN\dvtx X_n\ni u\}, \qquad u\in\bV,
\end{equation}
be the time that the node $u$ becomes an element of the BST sequence.
It is easy to see that the $\tau(u)$'s are finite
with probability 1.

%
\begin{theorem}\label{thm:DMconv}
Let $(X_n)_{n\in\bN}$ be the sequence of binary trees generated by the
BST algorithm with input a sequence
$(\eta_i)_{i\in\bN}$ of independent and identically distributed random
variables with $\cL(\eta_1)=\unif(0,1)$.

\begin{longlist}[(a)]
\item[(a)] With probability 1 the sequence $\unif(X_n)$ converges
weakly to a random probability measure $X_\infty$ on
$(\partial\bV,\cB(\partial\bV))$ as $n\to\infty$.

\item[(b)] For each $u\in\bV$, $u\neq\varnothing$, with $i:=\tau(u)-1$,
$\tau$ as in \eqref{eq:deftau},
and
\[
0=:\eta_{(i:0)} <\eta_{(i:1)}<\cdots<\eta_{(i:i)}<
\eta_{(i:i+1)}:=1
\]
the augmented order statistics associated with $\eta_1,\ldots,\eta_i$,
we have
\[
X_\infty(A_u) = \eta_{(i:j+1)} -\eta_{(i:j)}\qquad
\mbox{with } \eta_{(i:j)} <\eta_{i+1}<\eta_{(i:j+1)}.
\]

\item[(c)]
The random variables
\[
\xi_u:= \frac{X_\infty(A_{u0})}{X_\infty(A_u)},\qquad  u\in\bV,
\]
are independent, and $\cL(\xi_u)=\unif(0,1)$ for all $u\in\bV$.
\end{longlist}
\end{theorem}

\begin{pf}
Let $u$, $\tau(u)$, $i$ and $\eta_{(i:j)}$, $j=0,\ldots,i+1$, be as in
part (b) of the theorem.
The order property \eqref{eq:orderpropBST} of the labeled binary
search trees implies that for a node $v$ with label
$\eta_k$, $k>i$, the relation $v\ge u$ is equivalent to $\eta_{(i:j)}
<\eta_{k}<\eta_{(i:j+1)}$. Hence, by the law of
large numbers,
\begin{eqnarray*}
\lim_{n\to\infty} \unif(X_n) (\bar A_u) &=&
\lim_{n\to\infty} \frac{\#\{v\in X_n\dvtx v\ge u\}}{n}
\\
&=& \lim_{n\to\infty} \frac{\#\{i < k\le n\dvtx\eta_k\in(\eta
_{(i:j)}, \eta_{(i:j+1)})\}}{n}
\\
&=& \eta_{(i:j+1)} -\eta_{(i:j)}
\end{eqnarray*}
with probability 1 for every $u\in\bV$. In view of
\[
\{u\}= \bigl\{v\in\bar\bV\dvtx d(u,v)<2^{-|u|-1} \bigr\}\qquad \mbox{for all }
u \in\bV,
\]
the one-point sets with elements from $\bV$ are open in the topology on
$\bar\bV$. As $\unif(X_n)$ assigns at most the
value $1/n$ to such a set, it follows with the portmanteau theorem
\cite{Bill68}, page 11, that any limit point of this
sequence is concentrated on $\partial\bV$. Parts (a) and (b) of the
theorem now follow with the above general remarks
on weak convergence of probability measures on $(\bar\bV,\cB(\bar
\bV))$.

For the proof of (c) we use the following well-known fact: The
conditional distribution of $\eta_{i+1}$, given
$\eta_1,\ldots,\eta_i$ and given that the value lands in an interval
$I=(\eta_{(i:j)},\eta_{(i:j+1)})$ of the augmented
order statistics, is the uniform distribution on $I$, which implies
that $\unif(0,1)$ is the distribution of the normalized
distance $\xi_{u}$ to the left endpoint of $I$. For different $\eta
$-values these relative insertion positions
are independent, hence $\xi_u$, $u\in\bV$, are independent and
uniformly distributed on the unit interval.
\end{pf}

We note the following consequence of the representation in part (c) of
the theorem: For a fixed $u\in\bV$ let
\[
\varnothing=u(0)<u(1)<\cdots<u(k)=u
\]
with $|u(j)|=j$ for $j=0,\ldots,k$ be the path that connects $u$ to the
root node. We then have
%
\begin{eqnarray}
\label{eq:prodXinfty} X_\infty(A_u) = \prod
_{j=0}^{k-1}\tilde\xi_{u(j)}
\nonumber
\\[-8pt]
\\[-8pt]
\eqntext{\displaystyle\mbox{with }
\tilde\xi_{u(j)}:= %
\cases{ \xi_{u(j)}, &\quad $\mbox{if
}u(j+1)=u(j)0$, \vspace*{2pt}
\cr
1-\xi_{u(j)}, &\quad $\mbox{if
}u(j+1)=u(j)1$.}}
\end{eqnarray}
Note that the factors $\tilde\xi_{u(j)}$, $j=0,\ldots,k-1$, are
independent and that they all have
distribution $\unif(0,1)$.

Theorem \ref{thm:DMconv} confirms the view expressed in \cite{Woess2}, pages 191
and 218, that in specific cases embeddings
(or boundaries) can generally be obtained directly on using the then
available additional structure; here this turns out to
be the algorithmic representation of the Markov chain. However, there
are two additional benefits of the general theory:
First, because of the space--time property \eqref{eq:time-space} the
limit $X_\infty$ generates the tail $\sigma$-field
\[
\cT:= \bigcap_{n=1}^\infty\sigma \bigl(
\{X_m\dvtx m\ge n\} \bigr)
\]
associated with the sequence $(X_n)_{n\in\bN}$. This may serve as a
starting point for the unification of strong limit
theorems for functionals $(Y_n)_{n\in\bN}$, $Y_n=\Psi(X_n)$ of the
discrete structures: If $Y_n$ converges to $Y_\infty$
in a ``reasonable'' space, then the limit~$Y_\infty$, which is $\cT
$-measurable, must be a function of $X_\infty$; see, for example,
\cite{Kall}, Lemma 1.13. The second general result is extremely useful in
the context of the calculations that
arise in specific applications of the theory: The conditional
distribution of the chain $(X_n)_{n\in\bN}$ given the
value of $X_\infty$ is again a Markov chain, where the new transition
probabilities can be obtained from the limit value
and the old transition probabilities by a procedure that is known as
Doob's $h$-transform. In the present situation it
turns out that the conditional distribution of the BST chain, given
$X_\infty=\mu$, is the same as that of the DST chain
driven by $\mu$. We refer the reader to \cite{EGW} for details; the
last statement appears there only for a specific~$\mu$, but the generalization to an arbitrary probability measure $\mu$
in the boundary is straightforward. Roughly, the
embedded jump chains at the individual nodes are P\'olya urns; for
these the boundary has been obtained
in \cite{BK1964}, and from the general construction of the Doob--Martin
boundary it is clear that the outcome is
unaffected by the step from a Markov chain to its embedded jump chain.
We collect some consequences in the following
proposition, where
%
\begin{equation}
\label{eq:natFilt} \cF_n:= \sigma(X_1,\ldots,X_n),\qquad
n\in\bN,
\end{equation}
are the elements of the natural filtration of the BST chain.

%
\begin{proposition}\label{prop:conddistrnodes}
With the notation and assumptions as in Theorem \ref{thm:DMconv},
%
\begin{equation}
\label{eq:nodeBin} \cL \bigl(\sigma(X_n,u0) | \sigma(X_n,u)=k,
\xi_u=p \bigr) = \Bin(k-1,p) \qquad\mbox{if } k>0,
\end{equation}
and, for all $i,j\in\bN_0$,
%
\begin{equation}
\label{eq:nodeBeta} \cL \bigl(\xi_u | \sigma(X_n,u0)=i,
\sigma(X_n,u1)=j \bigr) = \Beta(i+1,j+1).
\end{equation}
Further, the variables $(\xi_u)_{u\in\bV}$ are conditionally
independent given $\cF_n$.
\end{proposition}

\section{Metric aspects}\label{sec:results}
All trees in this paper are subgraphs of the complete binary tree,
which has $\bV$ as its set of nodes and
$\{(\bar u,u)\dvtx u\neq \varnothing\}$ as its set of edges; in
particular, our trees are specified by their
node sets $x$. In a tree metric $d $ the distance of any two nodes
$u,v$ is the sum of the distances between
successive nodes on the unique path from $u$ to $v$, which means that
such a metric is given by
its values $d(\bar u,u)$, $u\in x$, $u \neq\varnothing$. For example,
the metric $d_\bV$ in
Section~\ref{subsec:DM} has $d_\bV(\bar u,u)=2^{-|u|-1}$, and the
canonical tree distance $d_{\mathrm{can}}$
is given by $d_{\mathrm{can}}(\bar u,u)=1$. For our trees the
prefix order further leads to
%
\begin{equation}
\label{eq:distviamin} d(u,v) = d(u,\varnothing) + d(v,\varnothing) - 2 d(u\wedge v,
\varnothing) \qquad\mbox{for all } u,v\in x.
\end{equation}
Metric trees may also be interpreted as graphs with edge weight, where
the edge $(\bar u,u)$ receives the weight
$d(\bar u,u)$.

Our aim in this section is to rephrase the convergence of the BST
sequence as a convergence of metric trees,
and to show that this view leads to convergence with respect to
stronger topologies. The situation here is
much simpler than for Aldous's continuum random tree where the
Gromov--Hausdorff convergence
of equivalence classes of metric trees is used; see \cite{EvansSF} and
the references given there. In fact,
the search trees considered here have node sets that grow monotonically
to the full $\bV$, so we may define
convergence of a sequence $((x_n,d_n))_{n\in\bN}$ of metric binary
trees to $(\bV,d_\infty)$ to mean that
%
\begin{equation}
\label{eq:treeconvweak} \lim_{n\to\infty} d_n(u,v) =
d_\infty(u,v) \qquad\mbox{for all } u,v\in\bV,
\end{equation}
which of course is equivalent to $ \lim_{n\to\infty} d_n(\bar u,u) =
d_\infty(\bar u,u)$ for all $u\in\bV$,
$u\neq\varnothing$. Note that $d_\bV$ and $d_{\mathrm{can}}$ are
both local metrics in the sense that $d(u,v)$ does
not depend on the tree $x$ as long as $u,v\in x$.

Motivated by the view in Section~\ref{subsec:DM} of finite and
infinite binary trees as probability measures $\mu$ on $(\bar\bV,\cB
(\bar\bV))$, we now introduce the (relative)
\emph{subtree size metric}, which assigns $\mu(\bar A_u)$ to the
distance of $\bar u$ and $u$, that is,
\[
d_x(\bar u,u) = \frac{\sigma(x,u)}{\sigma(x,\varnothing)} \qquad\mbox{for all } u\in x, u\neq
\varnothing,
\]
if $x\in\bB$, and
\[
d_\mu(\bar u,u) = \mu(A_u) \qquad\mbox{for all } u\in\bV, u
\neq\varnothing
\]
for the complete tree and a probability measure $\mu$ on $(\partial
\bV
,\cB(\partial\bV))$, where we assume that
$\mu(A_u)>0$ for all $u\in\bV$. Again, there is an algorithmic
motivation: In terms of the BST mechanism, the weight of
an edge $(\bar u,u)$ is the (relative) number of times this edge has
been traversed in the construction of the
tree. These metrics depend on their tree in a global manner.

With this terminology in place we may now rephrase the convergence in
Theorem \ref{thm:DMconv} as the
convergence in the sense of \eqref{eq:treeconvweak} of the finite
metric trees $(X_n,d_{X_n})$ to the infinite metric
tree $(\bV, d_{X_\infty})$, almost surely and as $n\to\infty$.

By construction the Doob--Martin compactification is the weakest
topology that allows for a continuous extension of the
functions $\bB\ni x\mapsto\sigma(x,u)/\break \sigma(x,\varnothing)$,
$u\in\bV$.
For the analysis of tree functionals stronger
modes of convergence turn out to be useful; for example, do we have
uniform convergence in \eqref{eq:treeconvweak}?
Also, subtree sizes decrease along paths leading away from the root
node, so we may consider a weight factor for the
distance of a node to its parent that depends on the depth of the node:
For all $\rho\ge1$, we define the
\emph{weighted subtree size metric} with \emph{weight parameter}
$\rho$ by
\[
d_{x,\rho}(\bar u,u):= \rho^{|u|}d_{x}(\bar u,u),\qquad
d_{\mu,\rho
}(\bar u,u):= \rho^{|u|}d_{\mu}(\bar u,u),
\]
in the finite and infinite case, respectively. Of course, with $\rho=1$
the subtree size metric reappears.

%
\begin{theorem}\label{thm:Hconv}
Let $\rho_0=1.26107\cdots$ be the smaller of the two roots of the
equation $ 2e\log(\rho)=\rho$, $\rho>0$.
Let $X_n$, $n\in\bN$ and $X_\infty$ be as in Theorem \ref{thm:DMconv}.
\begin{longlist}[(a)]
\item[(a)] For $\rho<\rho_0$, the metric space $(\bV,d_{X_\infty
,\rho
})$ is compact with probability 1.

\item[(b)] For $\rho> \rho_0$, the metric space $(\bV,d_{X_\infty
,\rho
})$ has infinite diameter with probability 1.

\item[(c)] For $\rho<\rho_0$, the metric spaces $(X_n,d_{X_n,\rho})$
converge uniformly to\break  $(\bV,d_{X_\infty,\rho})$
as $n\to\infty$ in the sense of
%
\begin{equation}
\label{eq:supconv}  \sup_{u,v\in X_n} \bigl| d_{X_n,\rho
}(u,v)-d_{X_\infty,\rho}(u,v)
\bigr| \to0 \qquad\mbox{ almost surely and in mean}.\hspace*{-35pt}
\end{equation}

\item[(d)] For $\rho> \rho_0$, and with $ d_{X_n,\rho}(\bar u,u):=0$
for $u\notin X_n$,
\[
\sup_{u,v\in\bV} \bigl| d_{X_n,\rho}(u,v)-d_{X_\infty,\rho}(u,v) \bigr| =
\infty\qquad\mbox{with probability 1}.
\]
\end{longlist}
\end{theorem}

\begin{pf}
We embed the metric trees into the linear space $\bL(0)$ of all
functions $f\dvtx\bV\setminus\{\varnothing\}\to\bR$ via
\[
x \mapsto f:= \bigl(u\mapsto d_x(\bar u,u) \bigr),\qquad  x\in\bB;
\]
probability measures $\mu$ on $(\bar\bV,\cB(\bar\bV))$ become elements
of $\bL(0)$ by identifying~$\mu$ with the function $u\mapsto\mu(A_u)$. In particular, we now
write $X_\infty(u)$ instead of $X_\infty(A_u)$.
For $\rho\ge1$ let $\bL(\rho)$ be the set of all $f\in\bL(0)$ with
\[
\|f\|_\rho:= \sum_{k=1}^\infty
\rho^k \max_{|u|=k} \bigl|f(u)\bigr| < \infty.
\]
Clearly, this gives a family of nested separable Banach spaces, with
\[
\bB\hookrightarrow\bL(\gamma)\subset\bL(\rho)\subset\bL(0) \qquad\mbox{for } 1\le
\rho<\gamma.
\]
We now show that, with the above identification,
%
\begin{eqnarray}
E \|X_\infty\|_\rho&< & \infty\qquad\mbox{if } \rho<
\rho_0, \label{eq:XinfinL1}
\\
P \Bigl( \sup_{u\in\bV} \rho^{|u|} X_\infty(u)=
\infty \Bigr) &= & 1 \qquad\mbox{if } \rho>\rho_0 \label{eq:XinfinL2}
\end{eqnarray}
and that, for $\rho<\rho_0$ and as $n\to\infty$,
%
\begin{equation}
\label{eq:liminL} \|X_n-X_\infty\|_\rho\to0 \qquad\mbox{almost surely and in mean.}
\end{equation}
Clearly, \eqref{eq:XinfinL1} implies that $X_\infty\in\bL(\rho)$ with
probability 1 if $\rho<\rho_0$.

The basis for our proof of \eqref{eq:XinfinL1} and \eqref{eq:XinfinL2}
is the connection of BST trees to branching
random walks, a connection that has previously been used by several
authors, especially for the analysis of the height
of search trees; see the survey \cite{DevGelbesBuch} and the references
given there. Let $u(k,j)$, $j=1,\ldots,2^k$, be
a numbering of the nodes from $\bV_k$ such that
\[
\beta \bigl(u(k,1) \bigr)<\beta \bigl(u(k,2) \bigr)<\cdots< \beta \bigl(u
\bigl(k,2^k \bigr) \bigr),
\]
with $\beta$ as defined in \eqref{eq:defbeta}. The key observation is
that the variables
\[
Y_{k,j}:= -\log X_\infty \bigl(u(k,j) \bigr), \qquad j=1,
\ldots,2^k,
\]
are the positions of the members of the $k$th generation in a branching
random walk
with offspring distribution $\delta_2$ and with
\[
Z:= \delta_{-\log\xi} +\delta_{-\log(1-\xi)},\qquad  \cL(\xi)=\unif(0,1)
\]
for the point process of the positions of the children relative to
their parent. Biggins~\cite{Biggins}
obtained several general results for such processes that we now
specialize to the present
offspring distribution and point process of relative positions. Let
\[
m(\theta):= E \biggl(\int e^{-\theta t} Z(dt) \biggr) = \frac
{2}{1+\theta}
\]
and
%
\begin{equation}
\label{eq:mudef} \tilde m(a):= \inf \bigl\{e^{\theta a} m(\theta )\dvtx\theta
\ge0 \bigr\} = 2ae^{1-a}.
\end{equation}
Note that
%
\begin{equation}
\label{eq:thetaopt} \tilde m(a) = m \bigl(\theta(a) \bigr)\qquad \mbox{with } \theta(a)=
\frac{1}{a}-1,
\end{equation}
and that, by definition of $\rho_0$,
%
\begin{equation}
\label{eq:rho2a} \rho<\rho_0\quad \Longleftrightarrow\quad \tilde m(\log\rho)<1.
\end{equation}
Finally, let $Z^{(k)}(t)$ be the number of particles in generation $k$
that are located to the left of $t$.

Now suppose that $\rho<\rho_0$. Let $\alpha:=(\rho+\rho_0)/2$ and
$\eta
:=\log(\alpha)$.
We adapt the upper bound argument in \cite{Biggins} to our present
needs: For all $\theta>0$
and $C > 1$, with $\gamma:=\log(C)$,
\begin{eqnarray*}
P \Bigl(\alpha^{k} \max_{|u|=k} X_\infty(u)>C
\Bigr) &=& P \Bigl( \min_{1\le j \le2^k} Y_{k,j} \le k\eta-\gamma
\Bigr)
\\
&\le& EZ^{(k)} \biggl(k \biggl(\eta-\frac{\gamma}{k} \biggr) \biggr)
\\
&\le&\exp \biggl(k \biggl(\eta-\frac{\gamma}{k} \biggr)\theta \biggr) m(
\theta)^k
\\
&=& C^{-\theta} \bigl(e^{\eta\theta}m(\theta) \bigr)^k.
\end{eqnarray*}
By \eqref{eq:rho2a}, $\tilde m(\eta)<1$. Choosing the optimal $\theta
=\theta(\eta)$, which
with \eqref{eq:thetaopt} is easily seen to be greater than $1$, leads to
\begin{eqnarray*}
E \Bigl(\alpha^{k} \max_{|u|=k} X_\infty(u)
\Bigr) &\le&1 + \int_1^\infty P \Bigl(
\alpha^{k} \max_{|u|=k} X_\infty(u)>x \Bigr)
\,dx
\\
&\le&1 + \tilde m(\eta)^k \int_1^\infty
x^{-\theta(\eta)} \,dx \le c,
\end{eqnarray*}
with a finite constant $c$ that does not depend on $k$. Hence
\[
\sum_{k=1}^\infty\rho^k E \Bigl(
\max_{|u|=k} X_\infty(u) \Bigr) \le c \sum
_{k=1}^\infty \biggl(\frac{\rho}{\alpha}
\biggr)^k < \infty,
\]
which in turn implies \eqref{eq:XinfinL1} by monotone convergence.

Suppose now that $\rho>\rho_0$, so that $\tilde m(\eta)>1$ by \eqref
{eq:rho2a} for $\eta:=\log\rho$.
By \cite{Biggins}, Theorem 2,
\[
\lim_{k\to\infty}\frac{1}{k} \log \bigl(\# \bigl\{1\le j
\le2^k\dvtx Y_{k,j}\le k \eta \bigr\} \bigr) = \log\tilde m(
\eta) > 0
\]
with probability 1. In particular, and again with probability 1,
\[
\exists k_0\ \forall k\ge k_0\ \exists u\in
\bV_k\qquad -\log X_\infty(u)\le k \log\rho.
\]
Clearly, this implies \eqref{eq:XinfinL2}.

For the proof of \eqref{eq:liminL} we first consider the random
variables $\sigma(X_n,u)$, $n\in\bN$, for some fixed
$u\in\bV$. We wish to relate these to $E[X_\infty(u)|\cF_n]$, with
$\cF
_n$ as in \eqref{eq:natFilt}. For this, we use
the representation of $X_\infty$ in terms of $(\xi_u)_{u\in\bV}$ given
in Section~\ref{subsec:DM}, together with
Proposition \ref{prop:conddistrnodes}. We may assume that $k:=|u|>0$.

The representation \eqref{eq:prodXinfty}, the conditional independence
of the $\tilde\xi$-variables given~$\cF_n$, and
the well-known formula for the first moment of beta distributions
together lead to
\[
E \bigl[X_\infty(u)|\cF_n \bigr] = \prod
_{j=0}^{k-1}E[\tilde\xi_{u(j)}|
\cF_n] = \prod_{j=0}^{k-1}
\frac{\sigma(X_n,u(j+1))+1}{\sigma(X_n,u(j)0) +
\sigma(X_n,u(j)1)+2}.
\]
In view of
\[
\sigma(x,u0)+ \sigma(x,u1)+1 = %
\cases{ \sigma(x,u), &\quad $\mbox{if } u
\in x$, \vspace*{2pt}
\cr
1, &\quad $\mbox{if } u\notin x$,} %
\]
the product telescopes to
%
\begin{equation}
\label{eq:telescprod2} E \bigl[X_\infty(u)|\cF_n \bigr] =
\frac{\sigma(X_n, u) +1}{n+1}\qquad \mbox{for all } u\in X_n.
\end{equation}
We now introduce
\[
Z_n\dvtx\bV\to\bR,\qquad  u\mapsto E \bigl[X_\infty(u)|
\cF_n \bigr].
\]
Then $(Z_n,\cF_n)_{n\in\bN}$ is a vector-valued martingale. For
$\rho
<\rho_0$ we have by part (a) of the theorem that
$X_\infty\in\bL(\rho)$ with probability 1 and that $E\|X_\infty\|
_\rho
<\infty$, hence $ Z_n\to X_\infty$ almost surely
and in mean in $\bL(\rho)$ by Proposition V-2-6 in \cite{NeveuMart}.

In our present representation of trees as functions on $\bV$ we have
\[
X_n(u) = %
\cases{ \displaystyle\frac{n+1}{n} Z_n(u)-
\frac{1}{n}, &\quad $\mbox{if } u\in X_n$,\vspace*{2pt}
\cr
0, &\quad $
\mbox{if } u\notin X_n$,} %
\]
which implies that $0\le X_n\le(1+n^{-1})Z_n$ for all $n\in\bN$. As
$X_n\to X_\infty$ pointwise with probability 1 by
Theorem \ref{thm:DMconv}
we can now use a suitable version of the dominated convergence theorem,
such as that given in \cite{Kall}, Theorem 1.21,
to obtain that $X_n$ converges to $X_\infty$ in $\bL(\rho)$ as $n\to
\infty$, again almost surely and in mean.

It remains to show that the tree statements in the theorem follow from
the linear space statements \eqref{eq:XinfinL1},
\eqref{eq:XinfinL2} and \eqref{eq:liminL}.

For (a) we prove that the limiting metric space is totally bounded.
From \eqref{eq:XinfinL1} and the definition of the
norm we obtain for any given $\varepsilon>0$ a $k=k(\varepsilon)\in
\bN$ such that
\[
\sum_{j=k}^\infty\rho^j \max
_{|u|=j} X_\infty(u) < \varepsilon,
\]
which by the definition of the weighted subtree size metric means that
all nodes $v$ with
$|v|\ge k$ have a distance from their predecessor at level $k$ that is
less than $\varepsilon$. As there are only
finitely many nodes of level less than $k$ this shows that the whole of
$\bar\bV$ may be covered by a finite
number of $\varepsilon$-balls. Of course, this argument is meant to be
applied to each element of a suitable set
of probability 1 separately.

For (b) we simply note that \eqref{eq:XinfinL2} implies that, with
probability 1,
\[
\sup_{u\in\bV,u\neq\varnothing} d_{X_\infty,\rho}(\bar u,u)=\infty
\]
if $\rho>\rho_0$. This also gives (d).\vadjust{\goodbreak}

Finally, for all $u\in\bV$, $u\neq\varnothing$,
\begin{eqnarray*}
\bigl| d_{X_n,\rho}(u,\varnothing)-d_{X_\infty,\rho}(u,\varnothing)\bigr | &\le&\sum
_{\varnothing\neq v\le u} \bigl| d_{X_n,\rho}(\bar v,v)-d_{X_\infty,\rho}(
\bar v,v) \bigr|
\\
&\le&\sum_{k=1}^{|u|} \rho^k \max
_{|v|=k} \bigl|X_n(v)-X_{\infty} (v)\bigr |
\\
&\le&\|X_n-X_\infty\|_\rho.
\end{eqnarray*}
The upper bound does not depend on $u$, hence (c) follows on
using \eqref{eq:distviamin}.
\end{pf}

We note that the convergence of metric trees considered in Theorem \ref
{thm:Hconv} implies the convergence with respect
to the Gromov--Hausdorff distance of the corresponding equivalence
classes of metric trees; see \cite{Burago}, Section~7.3.3.

%
\begin{figure}

\includegraphics{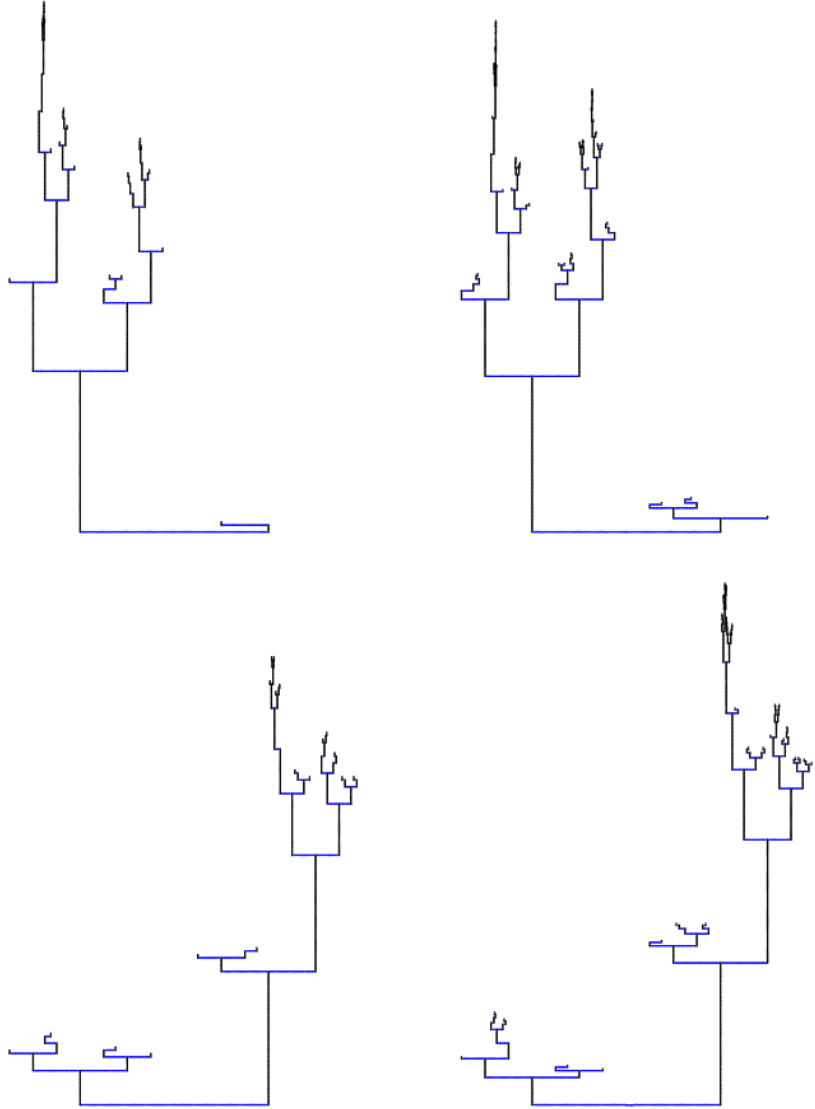}

\caption{The metric tree for the odd (upper part) and even (lower part)
$\pi$-data,
for $n=50$ (left) and $n=100$ (right), respectively; see text for
details.}\label{fig:pitrees}
\end{figure}

The subtree size metric also leads to a visualization of search trees:
We use the function $\beta$ defined
in \eqref{eq:defbeta} to map nodes to points in the unit interval, and
above the $x$-coordinate $\beta(u)$ we draw a
line parallel to the $y$-axis from $d_{X_n}(\bar u,\varnothing)$ to
$d_{X_n}(u,\varnothing)$. In order to obtain a visually
more pleasing result we may add lines that run parallel to the
$x$-axis, connecting nodes with the same parent. In
Figure~\ref{fig:pitrees} we have carried this out for the trees arising
from two separate input sequences for the BST
algorithm, with the data obtained from alternating blocks of length 10
of digits in the decimal expansion of
$\pi-3$. The upper part refers to the odd and the lower to the even
numbered blocks. In both cases we have given the trees
for $n=50$ and $n=100$, and with $\rho=1$. Vertically, the trees are
from the same distribution; moving horizontally to
the right, we have almost sure convergence.

\section{Tree functionals}\label{sec:application}
In this section we show how the above results can be used in connection
with the asymptotic analysis of tree functionals. Here is the recipe:
We start with a functional $Y_n=\Psi _n(X_n)$ of the trees, with
(deterministic) functions $\Psi_n$ on $\bB_n$ that have values in some
separable Banach space $(\bL,\|\cdot\|)$. We suspect that $Y_n$
converges almost surely to some limit variable $Y_\infty$ as $n\to
\infty $. We know that if this is the case, then
$Y_\infty=\Psi(X_\infty)$ for some $\Psi$ defined on $\partial\bB $ (as
always, almost surely). We do not know what $\Psi$ is, but if we manage
to rewrite the $\Psi_n$'s in terms of subtree sizes, then Theorem
\ref{thm:DMconv} may lead to an educated guess. On that basis we next
consider $\Phi_n(X_n)=E[\Psi(X_\infty)|\cF_n]$, assuming that $E\|\Psi
(X_\infty )\|<\infty$. This gives an $\bL$-valued martingale. By the
associated convergence theorem we then have that $\tilde
Y_n:=\Phi_n(X_n)$ converges to $Y_\infty$ almost surely and in mean.
Finally, a simple inspection of $\Phi _n-\Psi _n$ may reveal that
$\tilde Y_n-Y_n$ is asymptotically negligible---indeed, if $Y_n$
converges to $Y_\infty$, then $\tilde Y_n-Y_n$ \emph{must} tend to 0.

In the first three subsections we work out the details of the above
strategy for path lengths, for a tree index and for
an infinite dimensional tree functional. The final subsection is a
collection of remarks on other functionals and
related tree structures, indicating further applications of the method,
but also its limitations. The potential-theoretic
approach can provide additional insight; for example, we will relate a
martingale introduced in connection with tree
profiles to Doob's $h$-transform.

Throughout this section we abbreviate $X_\infty(A_u)$ to
$X_\infty(u)$.\newpage

\subsection{Path length}\label{subsec:PL}
The first tree functional we consider is the \emph{internal path length},
%
\begin{equation}
\label{eq:IPLa} \IPL(x):=\sum_{u\in x} |u|, \qquad x\in\bB,
\end{equation}
which may be rewritten as
%
\begin{equation}
\label{eq:IPLb} \IPL(x) =\sum_{u\in x,u\neq\varnothing} \sigma (x,u) =
\sum_{u\in x} \sigma(x,u) - \# x.
\end{equation}
Let
\[
H(0):=0,\qquad  H(n):=\sum_{i=1}^n
\frac{1}{i}\qquad \mbox{for all } n\in\bN,
\]
be the harmonic numbers. It is well known that
\[
\lim_{n\to\infty} \bigl( H(n) - \log n \bigr) =\gamma,
\]
where $\gamma\approx0.57722$ is Euler's constant.
We need two auxiliary statements; we omit the (easy) proofs.

%
\begin{lemma}\label{lem:intformula}
For all $i,j\in\bN_0$,
\[
\frac{\Gamma(i+j+2)}{\Gamma(i+1)\Gamma(j+1)}\int_0^1
x^i(1-x)^j \log(x) \,dx = H(i)-H(i+j+1).
\]
\end{lemma}

For a random variable
$\eta$ with distribution $\Beta(i+1,j+1)$ Lemma \ref{lem:intformula}
leads to
%
\begin{equation}
\label{eq:expxlogx} E \bigl(\eta\log(\eta) \bigr) = \frac
{i+1}{i+j+2} \bigl(H(i+1) -
H(i+j+2) \bigr).
\end{equation}
The next lemma is a summation by parts formula for binary trees.

%
\begin{lemma}\label{lem:telescope}
For any function $\psi:\bV\to\bR$,
\[
\sum_{u\in x} \bigl(\psi(u)-\psi(u0)-\psi(u1) \bigr) =
\psi(\varnothing) - \sum_{u\in\partial x} \psi(u) \qquad\mbox{for all
} x \in\bB.
\]
\end{lemma}

Major parts of the following theorem are known; we will give details
later in order to be able to
refer to the proof for a comparison of the methods used. Let
$(X_n)_{n\in\bN}$ be the BST chain, and
let $X_\infty$ be its limit, as in Theorem~\ref{thm:Hconv}.

%
\begin{theorem}\label{thm:IPL}
Let $C\dvtx(0,1)\to\bR$ be defined by
\[
C(s):= 1 + 2 \bigl( s\log(s) + (1-s)\log(1-s) \bigr).\vspace*{-9pt}
\]
\begin{longlist}[(a)]
\item[(a)] The limit
\[
Y_\infty:= \sum_{u\in\bV}X_\infty(u) C
\biggl(\frac{X_\infty
(u0)}{X_\infty(u)} \biggr)
\]
exists almost surely and in quadratic mean.

\item[(b)] As $n\to\infty$,
%
\begin{equation}
\label{eq:convIPL} \frac{1}{n} \IPL(X_n)-2\log n \to2\gamma- 4 +
Y_\infty,
\end{equation}
almost surely and in quadratic mean.
\end{longlist}
\end{theorem}

\begin{pf}
From the representation of $X_\infty$ given in Section~\ref{subsec:DM}, we know that the random variables
\[
\xi_u:=\frac{X_\infty(u0)}{X_\infty(u)},\qquad  u\in\bV,
\]
are independent and uniformly distributed on the unit interval, and
that $X_\infty(u)$ is a function of the $\xi_v$'s
with $v<u$. In particular, for all nodes $u$, the two factors in the
sum appearing in the definition of $Y_\infty$ are
independent. Let $\cG_k$ be the $\sigma$-field generated by the $\xi
_u$'s with $|u|\le k$, and put
\[
Y_k:= \sum_{u\in\bV, |u|\le k} X_\infty(u) C
\biggl(\frac
{X_\infty(u0)}{X_\infty(u)} \biggr),\qquad  k\in\bN.
\]
Then these properties lead to
\begin{eqnarray*}
E[Y_{k+1}| \cG_k] &=& Y_k + E \biggl[ \sum
_{|u|=k+1}X_\infty(u) C \biggl(\frac{X_\infty(u0)}{X_\infty(u)}
\biggr)\Big | \cG_k \biggr]
\\
&=& Y_k + \sum_{|u|= k+1} X_\infty(u)
EC(\xi_u)
\\
&=& Y_k,
\end{eqnarray*}
where we have used the fact that $EC(\xi_u)=0$. Further, with the same
arguments,
\begin{eqnarray*}
E \bigl[(Y_{k+1}-Y_k)^2| \cG_k
\bigr] &= &E \biggl[ \biggl( \sum_{|u|=k+1}
X_\infty(u) C(\xi_u) \biggr)^2 \Big|
\cG_k \biggr]
\\
&=& \sum_{|u|=k+1} X_\infty(u)^2 EC(
\xi_u)^2,
\end{eqnarray*}
so that
\[
E(Y_{k+1}-Y_k)^2 = \sum
_{|u|=k+1} EX_\infty(u)^2 EC(
\xi_u)^2.
\]
We also have $\kappa:=EC(\xi_u)^2<\infty$, and using \eqref
{eq:prodXinfty} we get
\[
EX_\infty(u)^2 = \bigl(E\xi_\varnothing^2
\bigr)^k = 3^{-k},
\]
so that
%
\begin{equation}
\label{eq:tailY} E(Y_{k+1}-Y_k)^2 =
2^k3^{-k} \kappa\qquad\mbox{for all } k\in\bN.
\end{equation}
Taken together these calculations show that $(Y_k,\cG_k)_{k\in\bN}$ is
an $L^2$-bounded martingale,
and an appeal to the corresponding martingale limit theorem completes
the proof of (a).
In particular, $Y_\infty$ is well defined, and even has finite second moment.

For the proof of (b) let $Z_n:=E[Y_\infty|\cF_n]$,
$n\in\bN$, so that $(Z_n,\cF_n)_{n\in\bN}$ is again a martingale
bounded in $L^2$. Our plan is to show that $Z_n$ is
sufficiently close to the transformed internal path length that appears
in \eqref{eq:convIPL}.

Using again the stochastic structure of $X_\infty$ we are thus led to
consider the conditional expectations
$E[X_\infty(u)|\cF_n]$ and $E[ C(\xi_u)|\cF_n]$, $u\in\bV$ and
$n\in\bN
$. From Proposition \ref{prop:conddistrnodes} we
know that, for all $u\in X_n$,
\[
\cL(\xi_u|\cF_n) = \Beta \bigl(\sigma(X_n,u0)+1,
\sigma(X_n,u1)+1 \bigr) ,
\]
and that the $\xi_u$'s are conditionally independent given $\cF_n$.
Hence Lemma \ref{lem:intformula}
can be applied [see also \eqref{eq:expxlogx}], resulting in
%
\begin{eqnarray}
\label{eq:condexpC} %
E \bigl[C(\xi_u)|\cF_n
\bigr] &=& 1 + \frac{2\tau(X_n,u0)+2\tau
(X_n,u1)}{\sigma(X_n,u0)+\sigma(X_n,u1)+2}
\nonumber
\\[-8pt]
\\[-8pt]
\nonumber
&&{}-2 H \bigl(\sigma(X_n,u0)+\sigma(X_n,u1)+2
\bigr),
\end{eqnarray}
where the function $\tau:\bB\times\bV\to\bR$ is given by
\[
\tau(x,u):= \bigl(\sigma(x,u)+1 \bigr) H \bigl(\sigma(x,u)+1 \bigr).
\]
For each fixed $n\in\bN$, almost sure convergence of $E[Y_k|\cF_n]$ to
$E[Y_\infty|\cF_n]$ as $k\to\infty$ follows
from
\[
\bigl\| E[Y_k|\cF_n]-E[Y_\infty|\cF_n]
\bigr\|_2 \le\| Y_k-Y_\infty \|_2,
\]
the upper bound in \eqref{eq:tailY} and the Borel--Cantelli lemma.
Together with the conditional independence of $X_\infty(u)$ and $C(\xi
_u)$ given $\cF_n$, this leads to
%
\begin{equation}
\label{eq:sumrepr} Z_n = \sum_{u\in\bV} E
\bigl[X_\infty(u) |\cF_n \bigr] E \bigl[C(\xi_u)|
\cF_n \bigr].
\end{equation}
From \eqref{eq:condexpC} we obtain $E[C(\xi_u)|\cF_n]=0$ for
$u\notin
X_n$, and, clearly,
%
\begin{equation}
\label{eq:splitnodes} \sigma(X_n,u0) + \sigma(X_n,u1) +1 =
\sigma(X_n,u)\qquad \mbox{for all } u\in X_n.
\end{equation}
Taken together, \eqref{eq:telescprod2}, \eqref{eq:condexpC}, \eqref
{eq:sumrepr} and \eqref{eq:splitnodes}
lead to
\[
Z_n = \sum_{u\in X_n}\frac{\sigma(X_n,u)+1}{n+1}
\biggl(1 + \frac{2\tau(X_n,u0)+2\tau(X_n,u1)}{\sigma(X_n,u)+1} -2 H \bigl(\sigma(X_n,u)+1 \bigr)
\biggr),
\]
which in turn gives
\[
Z_n = \frac{1}{n+1} \bigl(\IPL(X_n)+2n \bigr) -
\frac{2}{n+1} \sum_{u\in X_n} \bigl(
\tau(X_n,u)-\tau(X_n,u0)-\tau(X_n,u1) \bigr).
\]
Lemma \ref{lem:telescope} can be applied to the second sum, and the
assertion finally follows from
$\tau(X_n,\varnothing)=(n+1)H(n+1)$ and $\tau(X_n,u)=1$ for $u\in
\partial X_n$.
\end{pf}

Almost sure convergence of the standardized internal path length for
the BST sequence has been obtained
in \cite{RegnierQS}, and convergence in distribution, together with a
fixed point relation for the limit distribution,
in \cite{RoeQS}. Our method may been seen as an amalgamation of R\'
egnier's martingale approach and R\"osler's
approach, where the latter has come to be known as the contraction
method in the analysis of algorithms: We obtain
a strong limit, but we do not need to ``find the martingale'' (a task
familiar to many an applied probabilist).
The approach suggested in the present paper, to look at convergence of
the full objects via a suitable
completion of the state space of the underlying combinatorial Markov
chain, leads to a representation of the almost sure
limit. This gives the martingale by projection via conditional
expectations, and from the representation one can also read
off a fixed point relation for the distribution of the limit.

\subsection{The Wiener index}\label{subsec:WI}
The canonical graph distance $d_{\mathrm{can}} (u,v)$ of any two
nodes $u$ and $v$ in a finite connected graph $G$
with node set $V$ is the minimum length of a path (sequence of edges)
that connects $u$ and $v$ in $G$. The sum of these
distances is the \emph{Wiener index} of the graph,
%
\begin{equation}
\label{eq:WI} \WI(G):= \frac{1}{2}\sum_{(u,v)\in V\times V}
d_{\mathrm{can}} (u,v),
\end{equation}
introduced by the chemist H. Wiener. Some background together with
pointers to the literature is given in \cite{Nein02},
which is also our main reference in this subsection. Among other
results it is shown in \cite{Nein02} that for the BST
sequence $(X_n)_{n\in\bN}$ the rescaled Wiener indices,
\[
W_n:= \frac{1}{n^2} \WI(X_n) -2\log n,
\]
converge in distribution as $n\to\infty$.

Again, we project a suitable functional $\Psi(X_\infty)$ of the limit
tree $X_\infty$ to a function
$E[\Psi(X_\infty)|\cF_n]$ of $X_n$ that is sufficiently close to $W_n$.
This will give a strong limit theorem,
that is, it turns out that the rescaled Wiener indices in fact converge
almost surely for the random binary trees generated
by the BST algorithm for i.i.d. input, and it will also lead to a
representation of the limit $W_\infty$ as a function
of $X_\infty$.

We begin by rewriting the Wiener index in terms of subtree sizes,
similar to the transition from \eqref{eq:IPLa}
to \eqref{eq:IPLb} in the analysis of the internal path length. For a
binary tree $x$,
%
\begin{equation}
\label{eq:WIa} \sum_{(u,v)\in x\times x} |u\wedge v| = \sum
_{u\in x} \sigma(x,u)^2.
\end{equation}
This may be proved by induction, using the left and right subtrees in
the induction step;
see \cite{DennDiss}, page 70. Using \eqref{eq:distviamin}, \eqref
{eq:IPLb}, \eqref{eq:WI} and \eqref{eq:WIa}
we now obtain
%
\begin{equation}
\label{eq:WIrepr} \WI(X_n) = n \IPL(X_n) + n^2
- \sum_{u\in X_n} \sigma(X_n,u)^2.
\end{equation}
It is a benefit of working with
almost sure convergence that we can deal with the constituents on the
right-hand side of \eqref{eq:WIrepr}
separately (which means that we can make use of Theorem \ref{thm:IPL}),
whereas in connection with convergence
in distribution one needs to consider the joint distribution of $\IPL
(X_n)$ and $\WI(X_n)$; see \cite{Nein02}.

%
\begin{theorem}\label{thm:WI}
The series
%
\begin{equation}
\label{eq:defYinfty} Z_\infty:= \sum_{u\in\bV}X_\infty(u)^2
\end{equation}
converges almost surely and in quadratic mean, and, as $n\to\infty$,
%
\begin{equation}
\label{eq:thmWI} \frac{1}{n^2} \WI(X_n) -2\log n \to
W_\infty
\end{equation}
again almost surely and in quadratic mean, where the limit is given by
%
\begin{equation}
\label{eq:thmWIrepr} W_\infty:= 2\gamma-3 + Y_\infty-
Z_\infty,
\end{equation}
with $Y_\infty$ as in Theorem \ref{thm:IPL}.
\end{theorem}

\begin{pf}
Almost sure convergence in \eqref{eq:defYinfty} follows with
Theorem \ref{thm:Hconv}, and the moment calculations below
show that $EZ_\infty^2<\infty$. In particular,
\[
Z_n:= E[Z_\infty|\cF_n] \to Z_\infty
\]
almost surely and in quadratic mean. Again, the Markov property implies
that $Z_n$ can be written as a function of
$X_n$. In order to obtain this function we first consider a fixed node
$u\in\bV$.

From \eqref{eq:prodXinfty} we get
\[
X_\infty(u)^2 = \prod_{j=0}^{k-1}
\tilde\xi_{u(j)}^{ 2}.
\]
From \eqref{eq:nodeBeta} and the known formula for the second moment of
beta distributions we obtain,
considering the cases $u(j+1)=u(j)0$ and $u(j+1)=u(j)1$ separately,
\begin{eqnarray*}
&& E \bigl[\tilde\xi_{u(j)}^{ 2}|\cF_n \bigr]\\
&&\qquad =
\frac{(\sigma(X_n, u(j+1)) +1) (\sigma(X_n, u(j+1)) +2)} {
(\sigma(X_n, u(j)0) +\sigma(X_n, u(j)1) +2) (\sigma(X_n, u(j)0)
+\sigma(X_n, u(j)1) +3) }.
\end{eqnarray*}
Using the conditional independence statement in Proposition \ref
{prop:conddistrnodes}, we see that we have a
telescoping product again, so that
\[
E \bigl[X_\infty(u)^2|\cF_n \bigr] =
\frac{(\sigma(X_n, u) +1) (\sigma(X_n, u) +2)}{(n+1)(n+2)}
\qquad\mbox{for all } u\in X_n\cup\partial
X_n.
\]
The set $\bV\setminus X_n$ can be written as the disjoint union of the
subtrees rooted at the $n+1$ external nodes
of $X_n$, and we have
\[
E \bigl[\xi_u^2|\cF_n \bigr] = E \bigl[(1-
\xi_u)^2|\cF_n \bigr] = \tfrac{1}{3}\qquad
\mbox{for all } u\notin X_n.
\]
Therefore,
\begin{eqnarray*}
&&\sum_{u\notin X_n} E \bigl[X_\infty(u)^2|
\cF_n \bigr]
\\
&&\qquad= \frac{1}{(n+1)(n+2)}\sum_{u\in\partial X_n} \bigl(\sigma
(X_n,u)+1 \bigr) \bigl(\sigma(X_n,u)+2 \bigr) \sum
_{v\in\bV,v\ge u} \biggl(\frac{1}{3}
\biggr)^{|v|-|u|}
\\
&&\qquad= \frac{2}{(n+1)(n+2)}\sum_{u\in\partial X_n} \sum
_{k=0}^\infty2^k \biggl(\frac{1}{3}
\biggr)^k
\\
&&\qquad= \frac{6}{n+2}
\end{eqnarray*}
in view of $\sigma(X_n,u)=0$ for $u\in\partial X_n$. Taken together
this gives
\[
Z_n = \frac{1}{(n+1)(n+2)}\sum_{u\in X_n}
\bigl(\sigma(X_n,u)+1 \bigr) \bigl(\sigma(X_n,u)+2 \bigr)
+ \frac{6}{n+2}.
\]
Using \eqref{eq:IPLb} we get
\[
\sum_{u\in X_n} \bigl(\sigma(X_n,u)+1 \bigr)
\bigl(\sigma(X_n,u)+2 \bigr) = \sum_{u\in X_n}
\sigma(X_n,u)^2 + 3\cdot\IPL(X_n) + 5n
\]
so that, with \eqref{eq:WIrepr},
\[
\frac{1}{(n+1)(n+2)} \WI(X_n) = \frac{n}{(n+1)(n+2)} \IPL
(X_n) + \frac{n^2}{(n+1)(n+2)} - Z_n + R_n,
\]
where $R_n$ tends to 0 almost surely and in quadratic mean. From
Theorem \ref{thm:IPL} we know that
\[
\frac{1}{n} \IPL(X_n) - 2\log n \to2\gamma- 4 +
Y_\infty
\]
in the same sense. Combining the last two statements we obtain \eqref
{eq:thmWI}, with $W_\infty$ as in \eqref{eq:thmWIrepr}.
\end{pf}

\subsection{Metric silhouette}
In our third application we consider an infinite-dimensional tree functional.

Each element $v=(v_k)_{k\in\bN}$ of $\partial\bV$ defines a path
through a binary tree via the sequence
$(v(k))_{k\in\bN}$ of nodes given by $v(k)=(v_1,\ldots,v_k)$, $k\in
\bN
$. In \cite{GrSilh} the ``silhouette''
$\Sil(x)$ of $x\in\bB$ was introduced in an attempt to obtain a search
tree analogue of the famous Harris encoding of
simply generated trees: with each path $v$, we record its exit level
when passing through $x$, that is,
\[
\Sil(x) (v):= \min \bigl\{k\in\bN: v(k)\notin x \bigr\}, \qquad v\in \partial\bV.
\]
The tree silhouette can be visualized as a function on the unit
interval via the binary expansion
%
\begin{equation}
\label{eq:Phidef} \Phi\dvtx[0,1)\to\partial\bV, \qquad t\mapsto (v_k)_{k\in\bN}
\qquad\mbox{with } v_k:= \bigl\lceil2^{k+1}t \bigr\rceil- 2 \bigl
\lceil2^{k}t \bigr\rceil.
\end{equation}
It was shown in \cite{GrSilh} that for the BST chain $(X_n)_{n\in\bN}$
some smoothing is necessary to obtain an
interesting limit for the stochastic processes\break  $ (\Sil(X_n)(\Phi
(t)) )_{0\le t<1}$ as $n\to\infty$.

We have seen in the previous sections that for search trees it makes
sense to replace the canonical tree distance
implicit in the above definition of $\Sil(x)$ by the subtree size
metric. A corresponding variant
of the silhouette is the \emph{metric silhouette},
\[
\mSil(x) (v):= \sum_{k=1}^\infty\sigma
\bigl(x,v(k) \bigr),\qquad  v\in\partial\bV.
\]
Again, our aim is to obtain a strong limit theorem in the BST
situation, together with a representation of the limit as
a function of $X_\infty$. In addition, and going beyond the individual
arguments $v\in\partial\bV$, we regard
$\mSil(X_n)$ as a random function on $\partial\bV$. With $d_{\bV}$ as
in \eqref{eq:defdV} this is a compact and separable
metric space ($\bV$ is open in the completion $\bar\bV$ that we
introduced in Section~\ref{subsec:DM}). We write
$C(\partial\bV,d_{\bV})$ for the space of continuous functions
$f\dvtx\partial\bV\to\bR$. Together with
\[
\| f\|_\infty:= \sup_{v\in\partial\bV}\bigl |f(v)\bigr|,
\]
this is a separable Banach space.

Remember that the values of $X_\infty$ are probability measures on
$(\partial\bV,\cB(\partial\bV))$. Let
$\Sigma_\infty:\partial\bV\to[0,\infty]$ be defined by
\[
\Sigma_\infty(v):= -\int_{\partial\bV} \log_2
\bigl(d_{\bV}(u,v) \bigr) X_\infty(du),\qquad  v\in\partial\bV.
\]
This is the logarithmic potential of the random measure $X_\infty$ with
respect to $d_{\bV}$; see \cite{Woess1}, page 62.
Finally, we recall that a real function $f$ on the metric space
$(\partial\bV,d_{\bV})$ is said to be (globally) H\"older
continuous with exponent $\alpha$ if there exists a constant $C<\infty$
such that
\[
\bigl|f(u)-f(v)\bigr| \le C d_{\bV}(u,v)^\alpha\qquad\mbox{for all }u,v\in
\partial\bV.
\]

%
\begin{theorem}\label{thm:mSil}
Let $\alpha_0:= \log_2\rho_0=0.33464\ldots$ with $\rho_0$ as in
Theorem \ref{thm:Hconv}.
\begin{longlist}[(a)]
\item[(a)] $E \|\Sigma_\infty\|_\infty<\infty$.
\item[(b)] With probability 1, $\Sigma_\infty$ is H\"older continuous
with exponent $\alpha$ for all $\alpha<\alpha_0$.
\item[(c)] As $n\to\infty$,
\[
\biggl\| \frac{1}{n} \mSil(X_n) - \Sigma_\infty
\biggr\|_\infty\to0 \qquad\mbox{almost surely and in mean}.
\]
%
\end{longlist}
\end{theorem}

\begin{pf}
Because of $d_\bV(u,v)=2^{-|u\wedge v|}$ for all $u,v\in\partial\bV
$ we
have\break  $-\log_2 d_\bV(u,v)\in\bN_0$ and
\[
-\log_2 d_\bV(u,v) \ge k\quad \Longleftrightarrow\quad u\in
A_{v(k)}
\]
for all $k\in\bN$ so that
%
\begin{equation}
\label{eq:Sigma} \Sigma_\infty(v) = \sum_{k=1}^\infty
X_\infty \bigl(v(k) \bigr) \qquad\mbox{for all } v\in\partial\bV.
\end{equation}

Now let $\alpha$ be as in the statement of the theorem; we may assume
that \mbox{$\alpha>0$}. Let $\rho:=2^\alpha$. By
Theorem \ref{thm:Hconv} there exists a set of probability 1 such that
for all~$\omega$ in this set,
$C(\omega):=\|X_\infty(\omega)\|_\rho< \infty$. We fix such an
$\omega
$ and drop it from the notation. Because of
$X_\infty(u)\le C\rho^{|u|}$ for all $u\in\bV$ and \eqref
{eq:Sigma}, we
then have
$\Sigma_\infty(v)\le C\sum_{k=1}^\infty\rho^{-k}$ for all $v\in
\partial
\bV$,
which implies
\[
\|\Sigma_\infty\|_\infty\le\frac{1}{\rho-1} \|
X_\infty\|_\rho.
\]
In particular, $E\|\Sigma_\infty\|_\infty<\infty$ by \eqref
{eq:XinfinL1} in the proof of Theorem \ref{thm:Hconv}.

Similarly, if $u,v\in\partial\bV$ are such that $|u\wedge v|=k$, then
\begin{eqnarray*}
\bigl| \Sigma_\infty(u)-\Sigma_\infty(v)\bigr | &=& \sum
_{j=k+1}^\infty X_\infty \bigl(u(j) \bigr) + \sum
_{j=k+1}^\infty X_\infty \bigl(v(j)
\bigr)
\\
&\le&2 C \sum_{j=k+1}^\infty
\rho^{-j} = \frac{2C\rho^{-k}}{\rho-1} \le\frac{2C}{\rho-1}
d_{\bV}(u,v)^\alpha
\end{eqnarray*}
by definition of $\rho$. This proves (b).

For the proof of (c) we first consider the random functions $\Sigma_n$
defined by
\[
\Sigma_n(v):= E \bigl[\Sigma_\infty(v)|\cF_n
\bigr],\qquad v\in\partial\bV.
\]
With $J:=\Sil(X_n)(v)$ we get, using monotone convergence
for conditional expectations and $\cF_n$-measurability of $J$,
\begin{eqnarray*}
\Sigma_n(v) &=& \sum_{k=1}^J E
\bigl[X_\infty \bigl(v(k) \bigr)|\cF_n \bigr] + \sum
_{k=J+1}^\infty E \bigl[X_\infty \bigl(v(k)
\bigr)| \cF_n \bigr]
\\
&= &\sum_{k=1}^J \frac{\sigma(X_n,v(k))+1}{n+1} +
\frac{\sigma(X_n,v(J))+1}{n+1} \sum_{k=J+1}^\infty \biggl(
\frac{1}{2} \biggr)^{k-J}
\\
&= &\frac{1}{n+1} \bigl( \mSil(X_n) (v)+J \bigr) +
\frac{1}{n+1}.
\end{eqnarray*}
Here we have used our formula \eqref{eq:telescprod2} for $E[\Sigma
_\infty(u)|\cF_n]$ and its extension
to nodes outside $X_n$ that can be obtained as in the proof of
Theorem \ref{thm:WI}.

Let $h(x)=\max\{|u|: u\in x\}$ be the height of $x\in\bB$. Taking the
supremum over $v\in\partial\bV$ we get
\[
\biggl\| \Sigma_n - \frac{1}{n+1} \mSil(X_n)
\biggr\|_\infty\le\frac{h(X_n)+1}{n+1}.
\]
It is easy to show that the right-hand side converges to 0 with
probability 1 (see \cite{Dev86} for techniques and
results on the height),
so it remains to prove that $\Sigma_n$ converges almost surely and in
mean to $\Sigma_\infty$ in the separable
Banach space $ (C(\partial\bV,d_{\bV}),\|\cdot\|_\infty)$.
This, however, is again immediate
from the vector-valued martingale convergence theorem given in
\cite{NeveuMart}, page 104.
\end{pf}

%
\begin{figure}

\includegraphics{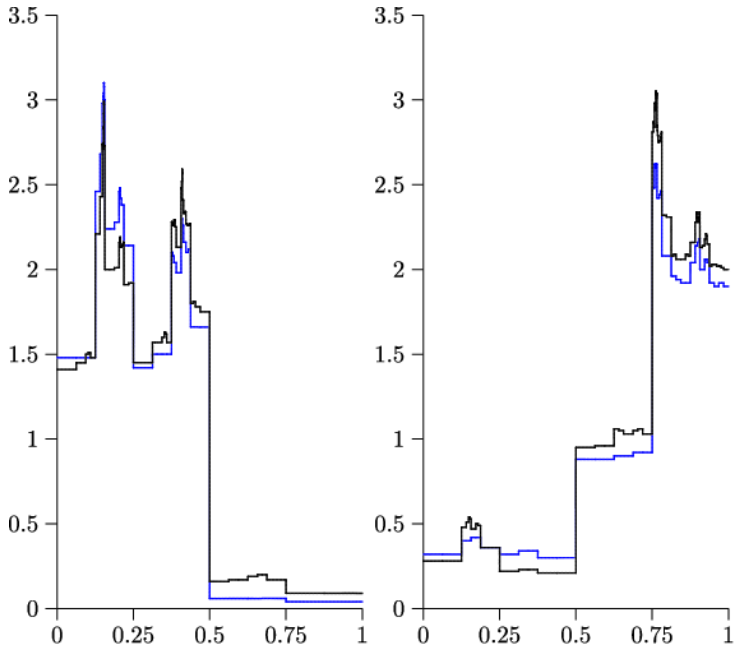}

\caption{The metric silhouette for the odd (left) and even (right)
$\pi$-data,
with $n=50$ (blue) and $n=100$ (black).}
\label{fig:metrSilh}
\end{figure}

Figure~\ref{fig:metrSilh} shows the metric silhouette for the trees in
Figure~\ref{fig:pitrees}. Note that the
continuity in Theorem \ref{thm:mSil} refers to the space $(\partial
\bV
,d_{\bV})$; for example,
$(t_n)_{n\in\bN}$ with $t_n=\frac{1}{2}+(-1)^n\frac{1}{n}$ for all $n
\in\bN$ is a Cauchy sequence with respect to
euclidean distance, but its inverse under the function $\Phi$ defined
in \eqref{eq:Phidef} that we used for the
illustration is not a Cauchy sequence in $(\partial\bV,d_{\bV})$.
Loosely speaking, the function~$\beta$ ``flattens''
the node set $\bV$.

\subsection{Other functionals and tree structures}\label{subsec:otherfunc}
The fill (or saturation) level $F(x)$ and height $H(x)$ of a tree $x\in
\bB$ are defined by
\[
F(x)=\max \bigl\{k\in\bN_0\dvtx\{0,1\}^k\subset x \bigr
\}, \qquad H(x)=\max\bigl\{|u|: u\in x \bigr\},
\]
respectively. For these tree functionals, the following asymptotic
results are well known:
%
\begin{equation}
\label{eq:fillheightconstants} \frac{F(X_n)}{\log n}\to\alpha_-,\qquad
 \frac{H(X_n)}{\log n}\to\alpha_+\qquad
\mbox{as } n\to\infty,
\end{equation}
both almost surely. Here $\alpha_-=0.373\ldots$ and $\alpha
_+=4.311\ldots$ are the two solutions of the equation
$x\log(2e/x)=1$. The survey \cite{DevGelbesBuch} gives details and
references, and explains the relation to branching
processes.

In situations such as these, where the almost sure limit is a constant,
projection on the sub-$\sigma$-fields
$\cF_n$ would simply return the constant, hence no simplification arises.

Both the fill level and height of a tree as well as its path length
(see Section~\ref{subsec:PL}) can be written as functionals
of the tree's node profile. Recall that $|u|$ denotes the length of
$u\in\bV$. Let
\[
v(x,k):= \# \bigl\{u\in\partial x\dvtx|u|=k \bigr\}, \qquad
w(x,k):= \#\bigl \{u\in x\dvtx|u|=k \bigr\}
\]
be the number of external (resp., internal) nodes of $x\in\bB$ at
depth $k$.
Applied to the BST sequence $(X_n)_{n\in\bN}$, this gives sequences
$(V_n)_{n\in\bN}$ and $(W_n)_{n\in\bN}$ of random
functions on the nonnegative integers via $V_n(k)=v(X_n,k)$ and
$W_n(k)=w(X_n,k)$, the external and internal node
profile of the binary search tree. Clearly,
\begin{eqnarray*}
F(X_n) &=& \min \bigl\{k\in\bN: V_n(k)>0 \bigr\} - 1,
\\
H(X_n) &=& \max \bigl\{k\in\bN: W_n(k)>0 \bigr\}
\end{eqnarray*}
and
\[
\IPL(X_n) = \sum_{k=1}^\infty
k W_n(k),
\]
so such profiles go some way toward a unifying approach to tree
functionals and indeed, they have been studied
extensively; see \cite{J-H2001,CDJ2001,CR2004,CKMR2005,FHN2006,DJN2008}.
A~crucial role
in \cite{J-H2001,CDJ2001,CR2004,CKMR2005} is played by a parametrized
family of martingales introduced in \cite{J-H2001}.
In order to connect this to the point of view of the present paper we
rephrase the basic idea using our terminology and
notation.

Fix some $z\in\bR_+$. The external profile $V_n$ of $X_n$ can be
regarded as the counting density of a
random finite measure on $\bN$ with total mass $n+1$ and value
\[
Y_n:= \sum_{k=1}^\infty
V_n(k) z^k = \sum_{u\in\partial
X_n}
z^{|u|}
\]
%
at $z$ of its generating function.
Let $v$ be the random node that is added to $X_n$ to obtain $X_{n+1}$.
It is easy to see that
\[
Y_{n+1} = Y_n + z^{|v|}(2z-1).
\]
In the BST mechanism the node $v$ is chosen uniformly at random from
the $n+1$ external nodes of $X_n$, hence
\begin{eqnarray*}
E[Y_{n+1}|\cF_n] &= & E \biggl[\sum
_{v\in\partial X_n} \bigl(Y_n+z^{|v|}(2z-1) \bigr)
1_{\{X_{n+1}=X_n\cup\{v\}\}} \Big|\cF_n \biggr]
\\
&=& Y_n E \biggl[\sum_{v\in\partial X_n}1_{\{X_{n+1}=X_n\cup\{v\}\}
}
\Big|\cF_n \biggr]
\\
&&{} + (2z-1) E \biggl[\sum_{v\in\partial X_n} z^{|v|}
1_{\{X_{n+1}=X_n\cup\{v\}\} }\Big |\cF_n \biggr]
\\
&= & Y_n + \frac{2z-1}{n+1}\sum_{v\in\partial X_n}
z^{|v|}
\\
&=& \frac{n+2z}{n+1} Y_n,
\end{eqnarray*}
which means that $(M_n,\cF_n)_{n\in\bN}$ with
\[
M_n:= C(n) Y_n, \qquad C(n):=\prod
_{k=1}^{n-1} \frac{k+1}{k+2z} \mbox{ for all } n\in
\bN,
\]
is a martingale. Obviously, the martingale is strictly positive
whenever $z>0$. Because of the space--time property
it can therefore be written as $M_n=h(X_n)$ with some positive harmonic
function $h$ on $\bB$, which depends on $z>0$,
and which in the present context is given by
\[
h(x) = C(\# x)\sum_{u\in\partial x}z^{|u|}.
\]
Moreover, the distribution $P^h$ of the corresponding $h$-transform,
which is the Markov chain with transition probabilities
\[
p^h(x,y) = \frac{1}{h(x)} p(x,y) h(y),\qquad  x,y\in\bB,
\]
is such that for all $n\in\bN$ the restriction $P^h_{\cF_n}$ of $P^h$
to $\cF_n$ has density
\[
\frac{dP^h_{\cF_n}}{dP_{\cF_n}} = \frac{1}{2z} h(X_n)
\]
with respect to the restriction $P_{\cF_n}$ to $\cF_n$ of the
distribution $P$ of the original BST chain.
Here we have used that both chains start with the tree $X_1=\{
\varnothing
\}$, and that $h(\{\varnothing\})=2z$.
A straightforward calculation yields
%
\begin{equation}
\label{eq:tiltedMC} p^h \bigl(x,x\cup\{v\} \bigr) = \frac{1}{n+2z}
\frac{\sum_{u\in\partial
x}z^{|u|} + z^{|v|}(2z-1)}{\sum_{u\in\partial x}z^{|u|}}
\end{equation}
for all $x\in\bB$, $v\in\partial x$.
Note that this agrees with the transition mechanism of the BST chain if
$z=1/2$ or $z=1$. For general $z>0$
a corresponding chain may be constructed by a marking mechanism that
makes use of an additional spine
variable. This idea was introduced in the context of branching
processes; for search trees it has been
used in \cite{CKMR2005}, to which paper we refer for more details. The
following direct construction of a
Markov chain with transitions as in \eqref{eq:tiltedMC} may be of
interest: Given $X_n$, we choose an external
node $u$ with probability proportional to~$z^{|u|}$. With probability
$(2z)/(n+2z)$ we then accept $u$ as the node
$v$ to be added to $X_n$; if $u$ is rejected, then $v$ is chosen
uniformly at random from the other $n$
external nodes of $X_n$.

In the first three subsections of the present section we began our
analysis by relating the functionals in question to
the subtree sizes. As the latter fully describe the tree this must also
be possible in the profile context.
For $x\in\bB$, $z>0$ let
\[
\Psi_z(x):= \sum_{u\in x} \sigma(x,u)
z^{|u|}.
\]
Each $v\in\bV_{k+1}$ with $\bar v\in x$ is either an internal or an
external node of $x$, which means that $v(x,k+1)=2
w(x,k) - w(x,k+1)$. Also, the number of internal nodes with depth at
least $k$ is the sum of all subtree sizes of the
nodes with level exactly equal to $k$, that is, $\sum_{u\in x,
|u|=k}\sigma(x,u)=\sum_{j=k}^\infty w(x,j)$. Taken together,
this gives
\[
\sum_{k=1}^\infty v(x,k) z^k =
\biggl(2z-3+\frac{1}{z} \biggr) \Psi_z(x) + \biggl(2-
\frac{1}{z} \biggr) \#x + 1,
\]
which leads to
\[
Y_n = \biggl(2z-3+\frac{1}{z} \biggr) \Psi_z(X_n)
+ \biggl(2-\frac{1}{z} \biggr) n + 1
\]
(note that the bracketed term vanishes for $z=1/2$ and $z=1$). This
could serve as the basis for an analysis along the
lines of the first three subsections. We do not pursue this here but
show instead that the general theory can be used
to obtain an interpretation of the function that represents the limit
of Jabbour's martingale in terms of the
Doob--Martin limit of the BST sequence: Recall that $M_n/(2z)$ is the
density associated with the change of measure from
$P_{\cF_n}$ to $P^h_{\cF_n}$. If the convergence $M_n\to M_\infty$ is
in $L^1$ (see below), then $M_\infty/(2z)$ is a
density of $P^h$ with respect to $P$. Thus we have $M_\infty=2z\Psi
(X_\infty)$, with $\Psi$ a density of the
distribution of $X_\infty$ under the transformed measure $P^h$ with
respect to the distribution of $X_\infty$ under the
original $P$.

It is shown in \cite{CKMR2005} that $L^1$-convergence holds if and only
if the parameter $z$ is inside a
specific bounded interval $I=(c_-,c_+)$, that $M_\infty\equiv0$ if
$z\notin I$, and that, with $\alpha_+,\alpha_-$ as
in \eqref{eq:fillheightconstants}, $c_-=\alpha_-/2$ and $c_+=\alpha
_+/2$. These two phase transitions are related to the
asymptotics of the maximum and minimum node size respectively at a
specific level of the limit $X_\infty$: If $z$ is too
small, then nodes close to the root are favored too much by $p^h$; if
$z$ is too large, then too much weight is given to nodes
far away from the root. In both cases $P^h$ is then singular with
respect to $P$. For the weighted subtree size
metric considered in Section~\ref{sec:results} only one of these
caveats matters in that node sizes must not be inflated
too much. Hence there is only one such phase transition, which should
be related to the height constant, and indeed, a
straightforward calculation shows that $\rho_0=(2e)/\alpha_+$.

Finally, let us mention that the approach toward strong asymptotics of
dynamic data structures that we have developed
in detail for binary search trees should be applicable in many related
situations. The necessary modifications may be
minor, such as for the discounted path length that appears in \cite
{GrSilh}, or straightforward, as for the random
recursive trees that are often treated in parallel with binary trees
(see, e.g., \cite{Nein02} for the Wiener index), or
they may be challenging, for example, when we wish to amplify the weak
convergence results for node depth profiles
obtained in \cite{DJN2008} for a wide class of trees to strong limit
theorems as we have done for the Wiener index in
Section~\ref{subsec:WI}. Of course, convergence in distribution and
convergence along paths are rather different
phenomena; see Figures~\ref{fig:pitrees} and \ref{fig:metrSilh}. It is
interesting that for a given dynamical structure
we may have a strong limit theorem (with nontrivial limit) for some
aspects (functionals), but not for others;
see \cite{DeGr3} for such results in connection with the subtree size
profile of binary search trees.

\section*{Acknowledgment}
I thank the referee for the stimulating comments,
and for
providing several valuable references.

%



\printaddresses

\end{document}